\documentclass[12pt]{article}
\usepackage{amssymb}
\usepackage{myart}

\normalbaselines
\input tcilatex
\begin{document}

\author{Yuri A Rylov}
\title{Crisis in the geometry development and its social consequences}
\date{Institute for Problems in Mechanics, Russian Academy of Sciences,\\
101-1, Vernadskii Ave., Moscow, 119526, Russia.\\
e-mail: rylov@ipmnet.ru\\
Web site: {$http://rsfq1.physics.sunysb.edu/\symbol{126}rylov/yrylov.htm$}\\
or mirror Web site: {$http://gasdyn-ipm.ipmnet.ru/\symbol{126}%
rylov/yrylov.htm$}}
\maketitle

\begin{abstract}
The reasons of the crisis in the contemporary (Riemannian) geometry are
discussed. The conventional method of the generalized geometries
construction, based on a use of the topology, leads to a overdetermination
of the Riemannian geometry. In other words, at the Riemannian geometry
construction one uses the needless information (topology), which disagrees
with other original axioms. The crisis manifests in the fact, that the
mathematical community cannot see and does not want to see the
overdetermination of the Riemannian geometry. The most geometers-topologists
deny the alternative method of the generalized geometry construction, which
does not uses the topology, because it does not contain theorems. Most
geometers see the geometry presentation as a set of definitions and
theorems. They cannot imagine the geometry presentation without customary
theorems. As a result the most clever topologists, which have acknowledged
the negligible role of the topology in the geometry construction and
inconsistency of the conventional method of the generalized geometry
construction, appear in the difficult situation (conflict with the
mathematical community).
\end{abstract}

\newpage 

\section{Introduction}

I wrote on crisis in geometry \cite{R2005}. The fact is that the
contemporary (Riemannian) geometry is overdetermined and hence it is
inconsistent. Although the character of inconsistency is known long ago, the
mathematical community as a whole does not want to acknowledge these
inconsistencies and bypasses them. The content of the crisis lies in the
fact that the mistakes in geometry (more exact in the construction of
generalized geometries) are not acknowledged and not corrected, but not in
the existence of mistakes in themselves. As a result the geometry is
developed in the direction, which leads in the blind alley. Why does such a
situation take place and which are consequences of the crisis? The presented
paper is devoted to this question.

Let us note that it the second crisis in the geometry. The first crisis took
place in the second half of the 19th century, when the mathematical
community did not want to acknowledge the reality of the non-Euclidean
geometries. This crisis had exhausted itself somehow, when the non-Euclidean
(Riemannian) geometry begun to be applied in the general relativity. But the
analysis of the crisis reasons has not been produced (at any rate, I know
nothing about such an analysis).

At first, I  connected an appearance of the crisis with the strongly crusted
preconception, that the straight is a one-dimensional line in any
generalized geometry. It followed from this statement that the
one-dimensional line is the most important object of any generalized
geometry. In its turn this fact leads to the conclusion that the topology,
founded on the concept of one-dimensional curve, lies in foundation of a
geometry, and it should be used in construction of a generalized geometry.
All conventional methods of the geometry construction are founded on the
essential employment of the topology. It means that at first one constructs
a topological space, and the geometry was constructed on the basis of the
topological space. Overcoming the preconception on the one-dimensionality of
the straight, one can  construct a more general and simpler method of the
generalized geometries construction. This alternative method is founded on
the deformation principle \cite{R2003}. This method admits, that the
straight may be a surface (tube), and it does not need the topology \cite%
{R2002} for construction of a generalized geometry. The substance of the
alternative method may be expressed by the words: "Any generalized geometry
may be obtained as a deformation of the proper Euclidean geometry. The
alternative method, applied for the Riemannian geometry construction, is
free of the well known defects of the conventional method, i.e. it does not
contain such properties of the Riemannian geometry as absence of the
absolute parallelism and the convexity problem. All scientists understood
that these properties of the Riemannian geometry were its defects, but they
could not overcome them without a use of the alternative method.

The fact, that the alternative method gives another results than the
conventional one, means that at least one of methods is erroneous. To use
the conventional method, one needs many different constraints, restricting a
possibility of its application. One needs the dimensionality of the space,
its continuity and a continuous coordinate system in it. For application of
the alternative method one does not need neither dimensionality, nor the
space continuity. The coordinate system is not needed also. One needs only
the distance function $\rho \left( P,Q\right) $, or the world function $%
\sigma \left( P,Q\right) =\frac{1}{2}\rho ^{2}\left( P,Q\right) $, which
satisfies the only condition 
\begin{equation}
\sigma \left( P,P\right) =0  \label{a1.1}
\end{equation}%
Even the symmetry condition 
\begin{equation}
\sigma \left( P,Q\right) =\sigma \left( Q,P\right)   \label{a1.2}
\end{equation}%
is not an obligatory condition \cite{R2002b}. The alternative method does
not contain any logical reasonings, whereas at application of the
conventional method of the generalized geometry construction one needs to
provide consistency of all used suppositions. It needs the significant
efforts. Consistency of all constraints may be provided not always. For the
Riemannian geometry the consistency is violated because of the supposition,
that the straight is an one-dimensional curve (or because of the application
of the topology at the geometry construction, what is essentially the same).
The fact is that at construction by means of the right (alternative) method
the geodesic, passing through the point $P$, in parallel to the vector at
the point $P$, is one-dimensional. However, in the case, when the straight
(geodesic) is passed through the point $P$, in parallel to the vector at the
other point $P_{1}$, it is not one-dimensional, in general \cite{R2006}. To
avoid the non-one-dimensionality of the straight, the absolute parallelism
is forbidden in the Riemannian geometry. Thus, it is evident that the
alternative method is true, whereas the conventional one is questionable. In
the case, when they lead to different results, one should to prefer the
alternative method.

I believed that the construction of generalized geometries by means of the
alternative method, founded on the deformation principle, admits one to
perceive the inadmissibility of the topology application at the generalized
geometry construction. But I appeared to be not right in the relation, that
I considered the one-dimensionality as an only preconception preventing from
the correct construction of the generalized geometry. It appears, that there
is once more preconception. I did not guess about this preconception,
because I did not meet it at the construction of the tubular geometry
(T-geometry), i.e. the geometry constructed on the basis of the deformation
principle. I met the preconception on the one-dimensionality of straight at
construction of the T-geometry and expended thirty years to overcome it \cite%
{R2005}. One can obtain a representation on these preconceptions from
consideration of following two syllogisms:

1. According to Euclid the straight is one-dimensional in the Euclidean
geometry, hence, the straight is one-dimensional in any generalized geometry.

2. Euclid constructed his geometry, formulating and proving theorems, hence,
any generalized geometry should be constructed, formulating and proving
theorems.

From the viewpoint of logic the considered syllogisms are not valid,
especially if one takes into account the opinion of the ancient Egyptians,
supposed that all revers flow northward. Their viewpoint may be considered
as a corollary of the third syllogism:

3. The great river Nile flows northward, hence all revers flow northward.

The structure of all three syllogisms is the same. They transform a partial
case into the general one without a sufficient foundations. Such a
generalization is connected with a narrowness of our experience. There is
only one (Euclidean) geometry, and the straight in this geometry is
one-dimensional. There is the only geometry, and it constructed by means of
formulation and proof of theorems. Finally, the ancient Egyptians knew only
one river Nile.

The logical inconsistency of syllogisms did not prevent mathematicians from
application of them at construction of any generalized geometry, because in
the given case the  mathematicians used in their practical activity the
associations, but not a logic. In all considered cases we deal with the only
object. In this case it is very difficult to distinguish, which properties
belong to the object in itself and which properties belong to the method of
the object description. When the method of description is considered to be a
property of the object in itself, we obtain a preconception. The first
syllogism is a reason of the preconception on one-dimensionality of the
straight. (The method of the straight description, used by Euclid, was
acknowledged as the property of the straight in itself). I knew about this
preconception, because I was forced to overcome it. The second syllogism
became to be a foundation of the preconception, that all activity of a
geometers consists in a formulation and a proof of theorems. (Euclid
described the geometry by means of theorems, hence description in terms of
theorems is a property of the geometry. If there are no theorems, then there
is no geometry.) I had met the second preconception at the following
circumstances.

I have submitted my report \cite{R2002} to the seminar of one of
geometric-topological chairs of the Moscow Lomonosov University. It took
place in 2001, and my paper \cite{R2002} was not yet published. I came to
the secretary of the seminar, told him about my wish to read a report and
gave him the text of the paper. Skimming my paper the secretary of the
seminar said thoughtfully: "How strange geometry! No theorems! Only
definitions! I believe that it will be not interesting for us." One of
leading geometers, which was forced duty bound to read my paper (the paper
was, submitted to the International conference on geometry, which took place
in Saint-Petersburg in the summer of 2002) said me in private conversation,
that he had understood nothing in my paper, because it contained neither
axioms, nor theorems. In that time I did not understand his objections. It
was quite unclear for me, how one cannot understand such a simple conception
as the T-geometry, which did not contain any logical reasonings. Only
essentially later I had understand, that here we deal with a preconception.
Mathematicians distinguish from other peoples in the relation, that they
perceive the geometry exclusively via its formalism. Further I shall try to
show formally, how theorems are transformed into definitions, and the demand
of the geometry representation in the form of a set of theorems is not more,
than a preconception.

Speaking about social consequences of the crisis, I keep in mind as follows.
Now the topology is considered to be the most perspective direction of the
geometry development. All best geometers are concentrated in the development
of this direction. Having overcame the preconception on the
one-dimensionality of the straight and discovering the alternative method of
the geometry construction, the topology appears to be a cul-de-sac direction
in the geometry development. Many papers on verification of the geometry on
the basis of topology appear to be depreciated. However, these papers are
considered to be important and perspective as long as the mathematicians do
not want to take into account the alternative method and do not acknowledge
inconsistencies in the Riemannian geometry.

Let us imagine the following situation. A young talented topologist N has
solved a very difficult topological problem. A prestige international prize
is declared for solution of this problem. Solving the problem and publishing
the results, the topologist N discovers suddenly, that the problem is set
incorrectly, because the topology is founded on concepts of the Riemannian
geometry, but the Riemannian geometry is overdetermined (i.e. it is
inconsistent). The topologist N is awarded this prize quite deservedly,
because one did not find mistakes in his solution. But, if the problem is
stated incorrectly, it may have no correct solution at all, because
different methods of solution may lead to different results. The
mathematical community considered that the prize has been awarded correctly,
but the topologist N disagree with the opinion of the mathematical
community, as far as he understand, that his work may be considered as
outstanding one, as long as the mathematical community did not discover the
weakness of the ground of the solved problem. What does the topologists N
do? If he is a conformist, he may affect ignorance, that he does not
discover any incorrectness in the statement of the problem and get quietly
the awarded deserved prize. (Whilom the inconsistency of the Riemannian
geometry will be discovered!). It is impossible to prove, that he has known
about inconsistency of the Riemannian geometry. However, the following
problem remain before the topologist N. Further some of his previous papers
will be depreciated. How and what for to work in the region of the
topological verification of the geometry, if this direction leads to the
blind alley?

However, if the topologist N is a scrupulous scientists, but not a
conformist, he has only one possibility: to keep away from the opinion of
the mathematical community and to abandon from the awarded prize. Should he
declare  the reason of his decision? It is not a simple question. Declaring 
the reason of his decision, the topologist N  conflicts with the
mathematical community, which does not acknowledge the Riemannian
inconsistency at this moment. In his time the great Gauss did not risk to
conflict with the mathematical community on the question of the
non-Euclidean geometry existence. (According to information of the Felix
Klein \cite{K37} one discovered unpublished manuscripts of Gauss on the
non-Euclidean geometry among his papers). However, Gauss worked in many
different branches mathematics, he may admit himself to neglect his works in
the non-Euclidean geometry. The contemporary topologists works, as a rule,
only in topology. Acknowledgment of inconsistencies in the Riemannian
geometry is a drama for a topologist. At the same time such an
acknowledgment is a very fearless act.

I suspect, that even if  the topologists N discussed inconsistency of the
Riemannian geometry with his colleagues, he may meet an incomprehension. I
have the following reason for such a suspicion. In the summer 2002 the
international conference, devoted to the anniversary of the known Russian
geometer A.D. Alexandrov, took place in Saint-Petersburg. Submitting the
corresponding report \cite{R2002b}, I arrived at the conference. However,
the main goal of my participation in the work of the conference was not my
report, which appeared to be interesting for nobody. I wanted to discuss
with leading geometers of Russia a possibility of the alternative method of
the Riemannian geometry construction, which concerns the geometry
foundation. It was important for me, because I am not a mathematician, but a
physicist - theorist. I succeeded to discuss this question with some leading
geometers, but nobody did not understood me (More exactly, only one person
had understood me, but he was not a geometer. He was simply a collaborator
of the Saint-Petersburg branch of Mathematical institute, where the
conference took place). The leading geometers told me very polite, that the
question of the geometry foundation had been solved many years ago by
Hilbert and other great geometers. In the present time the problem of the
geometry foundation is interesting for nobody. In this discussion only the
question on the alternative method of the geometry construction was
considered. The question of the Riemannian geometry inconsistency was not
discussed.

In other time I wanted to discuss the paper \cite{R2002} at the seminar of
the well known mathematician, who had usually a preliminary talk with the
potential speaker. During this talk I mentioned that the result of
application of the alternative method of the Riemannian geometry
construction disagree at some points with the results of application of the
conventional method. According to my representation, this fact must intrigue
the leader of the seminar in discussion of my report. However, I was told
immediately, that the Riemannian geometry cannot be inconsistent. Our talk
was finished after this declaration.

Summarizing, I should note that the more talented the topologist, the
earlier he faces the inconsistency of the Riemannian geometry, lying in the
foundation of the topology and makes the necessary conclusions. A failure in
his scientific career may be a result of this inconsistency. Let me note,
that the Riemannian geometry inconsistencies may be manifested only on the
sufficiently high level of the geometry (topology) development, when new
contradictions (other than the convexity problem and the fernparallelism
problem) appear. This does not concern the standard formalism of the
Riemannian geometry (metric tensor, covariant derivatives, curvature, etc.).

At first site the drama in the scientific career of the topologist N seem to
be unreal. However, it has taken place already, puzzling the mathematical
community, because nobody cannot imagine that the reason lies not in the
character of the topologist N, but in the crisis, which the talented
topologist has acknowledged earlier, than other scientists. His behavior in
the given situation was very dignified, although it was not clear for his
colleagues. I shall not denominate the name of this topologist, although
persons close to mathematics and topology may easily determine it.

Farther the pure mathematical questions will be considered, which explain my
statement on the overdetermination of the Riemannian geometry. The second
section is devoted to comparison of the conventional method of the
generalized geometry construction with the alternative one. In the third
section the question on the Riemannian geometry inconsistency is considered.
In the fourth section one considers why the mathematicians discard the
T-geometry, constructed on the basis of the deformation principle.

\section{Comparison of the conventional and alternative methods of the
geometry construction}

The conventional method of the Riemannian geometry is as follows. One
considers $m$-dimensional surfaces in the $n$-dimensional Euclidean space ($%
m<n$). Those properties of $m$-dimensional surfaces, which do not depend on
the dimension $n$ of the accommodating Euclidean space, declared to be an
internal geometry of the $m$-dimensional surface. It is the Riemannian
geometry of the $m$-dimensional space. The Riemannian geometries are
restricted by the constraint, that the Riemannian space can be embedded
isometrically into the Euclidean space of sufficiently large dimension. The
geometry which cannot be embedded isometrically in the Euclidean space,
cannot be constructed by this method. Besides the Riemannian geometries are
continuous, that is connected with the application of the continuous
coordinate system at the construction of the Riemannian geometry. A use of
the conventional method of the generalized geometry construction, i.e. the
geometry more general, than the Riemannian one, contains a series of
resricitions on the generalized geometry. In particular, such a constraint
is the embeddability of the space with the generalized geometry into the
Euclidean space of the sufficiently large dimension. Besides, such a
generalized geometry contains such a characteristic of the geometry as the
dimension, which is a natural number $n$. Necessity of some dimension $n$ at
the generalized geometry seems to be something as a matter-of-course,
although, in reality the dimension is a corollary of the applied method of
the generalized geometry construction, when one uses the concept of a
manifold, which is a coordinate system, consisting $n$ independent
coordinates. The fact, that the dimension is not a necessary property of the
generalized geometry (it is rather the means of description), follows from
the fact that there is an alternative method of the generalized geometry
description, where the concept of dimension may be not introduced.

The alternative method of the generalized geometry description is founded on
the deformation principle, which states that any generalized geometry can be
obtained as a result of a deformation of the proper Euclidean geometry. Any
deformation means a change of distance between the points of the space. Any
deformation of the proper Euclidean space generates some generalized
geometry. It is produced as follows. One proves the theorem, that the proper
Euclidean geometry may be formulated in terms and only in terms of the world
function \cite{R2002}. (The world function is a half of the squared distance
between the two points of the space). It follows from the theorem, that any
geometrical object $\mathcal{O}_{\mathrm{E}}$ and any statement $\mathcal{R}%
_{\mathrm{E}}$ of the proper Euclidean geometry $\mathcal{G}_{\mathrm{E}}$
can be expressed in terms of the world function $\sigma _{\mathrm{E}}$ of
the Euclidean geometry $\mathcal{G}_{\mathrm{E}}$ in the form $\mathcal{O}_{%
\mathrm{E}}\left( \sigma _{\mathrm{E}}\right) $ and $\mathcal{R}_{\mathrm{E}%
}\left( \sigma _{\mathrm{E}}\right) $ respectively. \textit{The set of all
geometrical objects $\mathcal{O}_{\mathrm{E}}\left( \sigma _{\mathrm{E}%
}\right) $ and relations $\mathcal{R}_{\mathrm{E}}\left( \sigma _{\mathrm{E}%
}\right) $ between them forms the Euclidean geometry $\mathcal{G}_{\mathrm{E}%
}$.} To obtain the corresponding relations of the generalized geometry $%
\mathcal{G}$ it is sufficient to replace the world function $\sigma _{%
\mathrm{E}}$ of the proper Euclidean geometry with the world function $%
\sigma $ of the generalized geometry $\mathcal{G}$: 
\[
\mathcal{O}_{\mathrm{E}}\left( \sigma _{\mathrm{E}}\right) \rightarrow 
\mathcal{O}_{\mathrm{E}}\left( \sigma \right) ,\qquad \mathcal{R}_{\mathrm{E}%
}\left( \sigma _{\mathrm{E}}\right) \rightarrow \mathcal{R}_{\mathrm{E}%
}\left( \sigma \right) 
\]%
\textit{Then the set of all geometrical objects $\mathcal{O}_{\mathrm{E}%
}\left( \sigma \right) $ and relations $\mathcal{R}_{\mathrm{E}}\left(
\sigma \right) $ between them forms the generalized geometry $\mathcal{G}$.}

The alternative method of the geometry construction supposes, that the
generalized geometry $\mathcal{G}$ is determined completely by its world
function, and any statement of the generalized geometry $\mathcal{G}$ can be
obtained from the corresponding statement of the Euclidean geometry. In this
case for construction of the generalized geometry $\mathcal{G}$ one does not
need a coordinate system. One does not need either its dimension, or any
another information on topology of the generalized geometry $\mathcal{G}$.
The dimension and topology (if they exist) can be obtained from the world
function $\sigma $ of the generalized geometry $\mathcal{G}$. The method of
determination of the dimension is a such one, that the dimension may be
different at different points of the space, or the dimension may be not
exist at all. It means that the definition of the topology and of the
dimension independently of the world function is inconsistent, in general.
Thus, the conventional method of the generalized geometry construction is
overdetermined. It contains too many axioms, which are not independent. The
suggested alternative method is insensitive to the continuity or
discreteness of the geometry, as far as it uses nowhere the limiting
transition or continuous coordinate system. A use of the conventional method
of the geometry construction, when one postulates some axiom system, which
determines the generalized geometry, appears to be ineffective, because it
is very difficult to provide compatibility of original axioms. The demand of
their compatibility imposes unnecessary constraints upon the obtained
generalized geometries.

For instance, a construction of the Riemannian geometry may be realized by
two methods. Using the conventional method on the basis of the metric
tensor, giving on the $n$-dimensional manifold, one obtains the Riemannian
geometry $\mathcal{G}_{\mathrm{R}}$. Using the alternative method based on
the deformation principle, one obtains $\sigma $-Riemannian geometry $%
\mathcal{G}_{\sigma \mathrm{R}}$. Prefix $\sigma $ means, that the
generalized geometry $\mathcal{G}_{\sigma \mathrm{R}}$ has the property of $%
\sigma $-immanence ($\sigma $-immanence is the property of the geometry to
be described completely by the world function $\sigma $). If the world
function $\sigma $ is the same in geometries $\mathcal{G}_{\mathrm{R}}$ and $%
\mathcal{G}_{\sigma \mathrm{R}}$, the obtained generalized geometries $%
\mathcal{G}_{\mathrm{R}}$ and $\mathcal{G}_{\sigma \mathrm{R}}$ are very
close, but they distinguish in some details. For instance, in the $\sigma $%
-Riemannian geometry $\mathcal{G}_{\sigma \mathrm{R}}$ there is the absolute
parallelism (fernparallelism), whereas in the Riemannian geometry $\mathcal{G%
}_{\mathrm{R}}$ one fails to introduce the absolute parallelism. The fact is
that the concept of parallelism of two vectors in the Riemannian geometry $%
\mathcal{G}_{\mathrm{R}}$ is transitive according to the geometry
construction \cite{R2002a,R2003}. In the $\sigma $-Riemannian geometry $%
\mathcal{G}_{\sigma \mathrm{R}}$ there is an absolute parallelism, which is
intransitive, but the parallelism of two vectors at the same point is
transitive. Thus, if one introduces the absolute parallelism in the
Riemannian geometry, it will be intransitive, in general, and incompatible
with the original statement on the transitivity of the parallelism in
Riemannian geometry. To avoid inconsistency, one declares that in the
framework of the Riemannian geometry one cannot introduce the absolute
parallelism of vectors, placed at different points of the space.

The considered example shows that the conventional method of the generalized
geometry construction is overdetermined. Furthermore, if the proper
Euclidean geometry is considered to be as a special case of the Riemannian
geometry, and one constructs it in the same way, as one constructs the
Riemannian geometry, the proper Euclidean geometry acquires an absurd
property. If some region of the space is nonconvex, the Riemannian geometry,
constructed in this region with the Euclidean metric tensor, is not the
Euclidean geometry, in general, because some distances in such a geometry
are determined not along the straight lines of the Euclidean geometry, but
along lines belonging partly to the boundary of the region. The obtained
Riemannian space cannot be embedded, in general, isometrically into the
Euclidean space, although the nonconvex region is a part of the Euclidean
space. This absurd result shows, that the system of axioms of the
conventional method of the generalized geometry construction is
overdetermined. An application of the overdetermined method of the
generalized geometry construction leads, in general, to inconsistency. The
form of these inconsistencies depends essentially on the method, how these
axioms are used.

\section{On inconsistency of the Riemannian geometry}

Some defects of the Riemannian geometry, were known years ago, but somehow
they are not perceived as inconsistency of the axiom system, which
determines the conventional construction of the Riemannian geometry.
(Apparently, it is conditioned with the absence of the alternative method of
the generalized geometry construction). For instance, if one constructs the
Riemannian geometry as a special case of the Riemannian one, using the
conventional method of the Riemannian geometry, the convexity problem
appears, consisting in the fact, that a nonconvex region of the Euclidean
plane cannot be embedded isometrically, in general, in the Euclidean plane,
from which it is cut. The result is evidently absurd. However, for the
convex region such an embedding is possible, mathematicians bypass this
defect, considering geometries only on convex manifolds. (For instance, the
book of A.D.Alexandrov "Internal geometry of convex surfaces"). Another
defect is the problem of the fernparallelism, i.e. absence of definition of
the parallelism of remote vectors in the Riemannian geometry. This fact is
not considered as a defect of the Riemannian geometry also. In reality these
defects are corollaries of the overdetermination of the Riemannian geometry,
i.e. at the Riemannian geometry construction one uses more axioms, than it
is necessary for the geometry construction. Some axioms are incompatible
between themselves, or they are compatible only at some constraints, imposed
on the geometry. In principle, this overdetermination may lead to another
contradictions, which are not known now. The overdetermination of the
conventional method of the Riemannian geometry construction has been
discovered after appearance of the alternative method, using essentially
less amount of information, which is necessary for the geometry
construction. The alternative method of the geometry construction does not
contain overdetermination, the convexity problem, the problem of
fernparallelism and other defects, which are corollaries of this
overdetermination.

The main preconception, preventing the contemporary geometry from its
development, is the statement that the straight is a one-dimensional line in
any generalized geometry (the straight has no width) has been traced to
Euclid. It is true in the Riemannian geometry for the straight, passing
through the point $P$, in parallel to a vector at the point $P$. However, in
the case, when the straight (geodesic) passes through the point $P$ in
parallel to a vector at another point $P_{1}$, it is not one-dimensional, in
general, \cite{R2006}. To avoid the non-one-dimensionality of the straight,
the fernparallelism (i.e. the concept of parallelism of two remote vectors)
is forbidden in the Riemannian geometry. The one-dimensionality of the
straight may not be included in axiomatics of a generalized geometry. The
one-dimensionality of the straight may not be used at the generalized
geometry construction, because it may appear to be incompatible with other
axioms. (Forbidding the fernparallelism, one did not think on a possible
non-one-dimensionality of the straight. Simply an ambiguity in the
definition of the parallelism contradicted to the axiomatics and seems to be
unacceptable). Forbidding fernparallelism and the geometry on non-convex
manifolds, one can avoid a manifestation of the Riemannian geometry
inconsistency. However, it does not avoid the inconsistency in itself,
because it may appear in other form. In the case of arbitrary generalized
geometry the character of the straight (one-dimensional line, or
multidimensional surface) is determined by the form of the world function,
and one does need to make suppositions on the one-dimensionality, or on
non-one-dimensionality of the straight. Furthermore, one cannot demand the
one-dimensionality of the straight, because it leads to an overdetermination
of conditions of the generalized geometry construction, and as a corollary
to their inconsistency, or to a restriction of the class of possible
geometries. Application of the topology at the generalized geometry
construction supposes that a curve (and its kind a straight) is
one-dimensional. It is a reason, why the topology may not be used at the
generalized geometry construction. Overcoming of this preconception was very
difficult. I personally needed almost thirty years for overcoming this
preconception (a description of the path of this overcoming can be found in 
\cite{R2005}), although I have an experience of successful overcoming of
like preconception in other branches of physics.

Although the overcoming of the preconception on the one-dimensionality of
the straight was difficult, I hoped, that, being elucidated it will be
perceived and overcame by the mathematical community. It appeared, that I
was wrong in this relation. The overcoming of the perception was very
difficult ( by the way, it was difficult also in the case of the
preconception on the statistical description). Apparently, overcoming of
preconceptions is always difficult. Besides, it appeared, that besides of
the main preconception on the one-dimensionality of the straight there are
another preconceptions. For instance, one supposes, that any presentation of
a geometry is a set of axioms and theorems and that the geometer's activity
is a formulation and a proof of theorems.

\section{Why mathematicians do not acknowledge \newline
T-geometry constructed on the basis of the \newline
deformation principle?}

The relation of the most mathematicians to the alternative method of the
geometry construction, which is based on the deformation principle is
negative as a rule, although they cannot contradict anything against this
method. Indeed, it as very difficult to contradict anything, because of the
simplicity of the alternative method, which does not contain any
suppositions except for the sufficiently evident deformation principle.
Nevertheless, the refereed mathematical journals reject papers on the
generalized geometry construction founded on basis of the deformation
principle. In the paper \cite{R2005} it is described in detail, how such a
rejection is produced. My attempts of reporting these papers at the seminars
of topologists - geometers are rejected also. However, my papers were
reported and discussed at the seminar on the geometry as a whole in the
Moscow Lomonosov University. This seminar was founded by N.V. Efimov.

In the present time the topological approach to the problem of the
generalized geometry construction dominates. One supposes, that it is
necessary at first to construct a proper topological space. In the
topological space one introduces a metric and constructs a generalized
geometry. Construction of the generalized geometries on the basis of the
deformation principle depreciates essentially papers on geometry, and the
mathematicians, constructing generalized geometries by the conventional
method, were not enthused over the prospect. As far as they cannot suggest
any essential objection against the deformation principle because of its
simplicity and effectiveness, they use another methods of resistance,
generating doubts in the scientific scrupulosity of advocates of the
conventional method. One of them is described in \cite{R2005} in detail.

Of course, I have known, how strongly we believe, that the straight line is
always an one-dimensional geometrical object, because this preconception
retarded my discovery of the T-geometry almost with thirty years. However,
discovery of the non-one-dimensionality of the straight and incomprehension
of such a possibility, when the non-one-dimensionality has been already
discovered, seemed to be different things for me. It was appeared that a
perception of the deformation principle is the more difficult, the better
the person knows the geometry in its contemporary formal presentation. The
fact is that the mathematicians perceive all things via the formalism. The
associative perception like a reference to the deformation principle conveys
little to them.

The strangers and the mathematicians perceive the geometry presentation as a
formulation and proof of different geometric theorems. Replacement of
theorems with definitions seems to them as anything quite obscure. However,
in reality, the formulation and the proof of theorems is only one of the
methods of work with different geometric statements (axioms and theorems).
The proof of theorems is the most laborious part of work. As a result one
gets an impression, that the geometry presentation is a proof of theorems
(if there are no theorems, there is no geometry). But there are another
methods of work with the geometric statements, which does not distinguish
between the axioms and theorems.(There is no necessity to repeat the work of
Euclid, one should use results of his work.) At such a method of work the
theorems are replaced by definitions, and a necessity of their proof falls
off.

I shall explain this in the example of the cosine theorem, which states%
\begin{eqnarray}
\left\vert \mathbf{BC}\right\vert ^{2} &=&\left\vert \mathbf{AB}\right\vert
^{2}+\left\vert \mathbf{AC}\right\vert ^{2}-2\left( \mathbf{AB.AC}\right) 
\nonumber \\
&=&\left\vert \mathbf{AB}\right\vert ^{2}+\left\vert \mathbf{AC}\right\vert
^{2}-2\left\vert \mathbf{AB}\right\vert \left\vert \mathbf{AC}\right\vert
\cos \alpha  \label{a3.1b}
\end{eqnarray}%
where the points $A,B,C$ are vertices of a triangle, $\left\vert \mathbf{BC}%
\right\vert $, $\left\vert \mathbf{AB}\right\vert $, $\left\vert \mathbf{AC}%
\right\vert $ are lengths of the triangle sides and $\alpha $ is the angle $%
\measuredangle BAC$. The relation (\ref{a3.1b}) is the cosine theorem which
is proved on the basis of the axioms of the proper Euclidean geometry.

Using expression of the length of the triangle side $\mathbf{AB}$ via the
world function $\sigma $%
\begin{equation}
\left\vert \mathbf{AB}\right\vert =\sqrt{2\sigma \left( A,B\right) }
\label{a3.2}
\end{equation}%
we may rewrite the relation (\ref{a3.1b}) in the form%
\begin{equation}
\left( \mathbf{AB.AC}\right) =\sigma \left( A,B\right) +\sigma \left(
A,C\right) -\sigma \left( B,C\right)  \label{a3.3}
\end{equation}%
The relation (\ref{a3.3}) is a definition of the scalar product $\left( 
\mathbf{AB.AC}\right) $ of two vectors $\mathbf{AB}$\textbf{\ }and $\mathbf{%
AC}$ in the T-geometry. Thus, the theorem is replaced by the definition of a
new concept (the scalar product), which is not connected directly with the
concept of the linear space.

Another example the Pythagorean theorem for the rectangular triangle $ABC$
with the right angle $\measuredangle BAC$, which is written in the form 
\begin{equation}
\left\vert \mathbf{BC}\right\vert ^{2}=\left\vert \mathbf{AB}\right\vert
^{2}+\left\vert \mathbf{AC}\right\vert ^{2}  \label{a3.4}
\end{equation}%
In T-geometry instead of the theorem (\ref{a3.4}) we have a definition of
the right angle$\measuredangle BAC$. In terms of the world function this
definition has the form. The angle $\measuredangle BAC$ is right, if the
relation%
\begin{equation}
\sigma \left( A,B\right) +\sigma \left( A,C\right) -\sigma \left( B,C\right)
=0  \label{a3.5}
\end{equation}%
takes place

Thus, we see that theorems of the proper Euclidean geometry are replaced by
definitions of T-geometry.

From the formal viewpoint the difference between the conventional method of
description and the alternative method may be described as follows. The
conventional method of the proper Euclidean geometry may be described as a
set $\mathcal{S}_{\mathrm{E}}\left( \mathcal{A}\left[ \mathcal{R}_{\mathrm{E}%
}\right] \right) $ of algorithms $\mathcal{A}\left[ \mathcal{R}_{\mathrm{E}}%
\right] $, acting on operands $\mathcal{R}_{\mathrm{E}}$, where operands $%
\mathcal{R}_{\mathrm{E}}$ are geometrical objects or relations of the
Euclidean geometry. The operands $\mathcal{R}_{\mathrm{E}}$ depend on
parameters $\mathcal{P}^{n}\equiv \left\{ P_{0},P_{0},..,P_{n}\right\} $,
where $P_{0},P_{0},..,P_{n}$ are points of the space. 
\begin{equation}
\mathcal{R}_{\mathrm{E}}=\mathcal{R}_{\mathrm{E}}\left( \mathcal{P}%
^{n}\right)  \label{a3.6}
\end{equation}%
As far as all geometrical objects and relations $\mathcal{R}_{\mathrm{E}}$
can be expressed via the world function $\sigma _{\mathrm{E}}$ of the
Euclidean geometry $\mathcal{G}_{\mathrm{E}}$, one may rewrite the relation (%
\ref{a3.6}) in the form 
\begin{equation}
\mathcal{R}_{\mathrm{E}}=\mathcal{R}_{\mathrm{E}}\left( \mathcal{P}%
^{n}\right) =\tilde{\mathcal{R}}_{\mathrm{E}}\left[ \mathcal{\sigma }_{%
\mathrm{E}}\left( \mathcal{P}^{n}\right) \right]  \label{a3.7}
\end{equation}%
where $\mathcal{\sigma }_{\mathrm{E}}\left( \mathcal{P}^{n}\right) $ is a
set of world functions $\sigma _{\mathrm{E}}\left( P_{i},P_{k}\right) ,$ $%
P_{i},P_{k}\in \mathcal{P}^{n}$.

Taking into account the relation (\ref{a3.7}), one can represent the set $%
\mathcal{S}_{\mathrm{E}}\left( \mathcal{A}\left[ \mathcal{R}_{\mathrm{E}}%
\right] \right) $ of all algorithms $\mathcal{A}\left[ \mathcal{R}_{\mathrm{E%
}}\left( \mathcal{P}^{n}\right) \right] $ in the form of the set $\mathcal{S}%
_{\mathrm{E}}\left( \mathcal{A}\left[ \mathcal{R}_{\mathrm{E}}\left( 
\mathcal{P}^{n}\right) \right] \right) =\mathcal{S}_{\mathrm{E}}\left( 
\mathcal{A}\left[ \tilde{\mathcal{R}}_{\mathrm{E}}\left[ \mathcal{\sigma }_{%
\mathrm{E}}\left( \mathcal{P}^{n}\right) \right] \right] \right) =\mathcal{S}%
_{\mathrm{E}}\left( \mathcal{A}\left[ \tilde{\mathcal{R}}_{\mathrm{E,}%
\mathcal{P}^{n}}^{\prime }\left[ \mathcal{\sigma }_{\mathrm{E}}\right] %
\right] \right) $ of all algorithms $\widetilde{\mathcal{A}}\left[ \mathcal{%
\sigma }_{\mathrm{E}}\right] =\mathcal{A}\left[ \tilde{\mathcal{R}}_{\mathrm{%
E,}\mathcal{P}^{n}}^{\prime }\left[ \mathcal{\sigma }_{\mathrm{E}}\right] %
\right] $. In the Euclidean geometry the set of all algorithms $\mathcal{S}_{%
\mathrm{E}}\left( \mathcal{A}\left[ \mathcal{R}_{\mathrm{E}}\left( \mathcal{P%
}^{n}\right) \right] \right) $ takes the form $\mathcal{S}_{\mathrm{E}%
}\left( \widetilde{\mathcal{A}}\left[ \tilde{\mathcal{R}}_{\mathrm{E,}%
\mathcal{P}^{n}}^{\prime }\left[ \mathcal{\sigma }_{\mathrm{E}}\right] %
\right] \right) $. In the generalized geometry $\mathcal{G}$, described by
the world function $\mathcal{\sigma }$, the set of all algorithms takes the
form $\mathcal{S}_{\mathrm{E}}\left( \widetilde{\mathcal{A}}\left[ \tilde{%
\mathcal{R}}_{\mathrm{E,}\mathcal{P}^{n}}^{\prime }\left[ \mathcal{\sigma }%
\right] \right] \right) $. It means that the set of all algorithms $\mathcal{%
S}_{\mathrm{E}}\left( \widetilde{\mathcal{A}}\left[ \tilde{\mathcal{R}}_{%
\mathrm{E,}\mathcal{P}^{n}}^{\prime }\left[ \mathcal{\sigma }\right] \right]
\right) $ for any generalized geometry is obtained from the set $\mathcal{S}%
_{\mathrm{E}}\left( \widetilde{\mathcal{A}}\left[ \tilde{\mathcal{R}}_{%
\mathrm{E,}\mathcal{P}^{n}}^{\prime }\left[ \mathcal{\sigma }_{\mathrm{E}}%
\right] \right] \right) $ of all algorithms for the Euclidean geometry by
means of replacement of the operand of algorithms $\sigma _{\mathrm{E}%
}\rightarrow \sigma $. The form of all algorithms is not changed. This is
the formal description of the action of the deformation principle. For
construction of the generalized geometry one needs no theorems and no
logical reasonings. The main part of the work of construction of the
generalized geometry consists in transformation of known Euclidean
algorithms $\mathcal{A}\left[ \mathcal{R}_{\mathrm{E}}\left( \mathcal{P}%
^{n}\right) \right] $ to the form $\widetilde{\mathcal{A}}\left[ \tilde{%
\mathcal{R}}_{\mathrm{E,}\mathcal{P}^{n}}^{\prime }\left[ \mathcal{\sigma }_{%
\mathrm{E}}\right] \right] $, where all algorithms are represented in terms
of the Euclidean world function.

In the T-geometry a construction of new generalized geometry is produced by
the same algorithms, as in the proper Euclidean geometry. One replaces only
the world function, i.e. operand of algorithms $\sigma _{\mathrm{E}%
}\rightarrow \sigma $. One does not need theorems which connect different
objects of geometry. One uses the definition of the geometrical object or the%
\textrm{\ }relation of the Euclidean geometry, expressed directly via world
function. For instance, the scalar product of two vectors is defined by the
relation (\ref{a3.3}) in all T-geometries. The right angle $\measuredangle
BAC$ is defined by the relation (\ref{a3.5}) in all T-geometries. There are
no necessity to introduce couplings (theorems) between different geometrical
concepts, as far as they are expressed via world function.

In T-geometry the main problem is an obtaining of expressions for different
geometrical concepts and objects via the world function. As far as these
expressions are the same for all T-geometries, it is sufficient to obtain
these expressions in the framework of the proper Euclidean geometry. Thus,
instead of formulation of different theorems for different generalized
geometries, the geometer must express all concepts and objects of the proper
Euclidean geometry via the world function.

The formalism of T-geometry distinguishes essentially from the conventional
formalism of the generalized geometry construction by the object of
consideration. In the conventional formalism the geometrical objects and
primarily the straight are objects of consideration. The conventional
construction of the generalized geometry is a repetition of the Euclid's
construction, which is produced on the basis of other axioms. The true
choice of axioms is the main problem of the conventional method of the
geometry construction. In the alternative method, based on the deformation
principle, the object of consideration is the world function (but not
geometrical objects), i.e. a function of two points of the space. It is
essentially simpler object of consideration. However, the rules of work
(logic of investigation \cite{R2005a}) with the world function are very
complicated. It is very complicated to guess these rules, because they are
determined by the Euclidean geometry. Primarily one guesses the rules of
work with the most simple object: the straight line. As a matter of fact,
they were not guessed, they have been taken from the Euclidean geometry.
But, to take them, the Euclidean geometry is to be presented in terms of the
world function. It means that it was necessary to construct the formalism on
the basis of the world function. Construction of the formalism, founded on
the world function, was difficult. It needed about forty years. Thirty years
of these forty years were lost for the overcoming the preconception on the
one-dimensionality of straight. Different stages of this formalism
construction are described in \cite{R2005}.

Apparently, the most difficult and uncustomary for geometers - professionals
are the facts, that a new object of consideration (world function) appears
and the conventional objects of consideration (geometrical objects) turn
into the logic (algorithm) of investigation. This transition from one method
of investigation to another is difficult for perception of the geometers -
professionals. However, it is essentially easier for non- professionals, who
do not know, or know skin-deep the conventional method of the geometry
construction. They do not need to compare both methods, overcoming complex
of concepts of the conventional method and the formalism, connected with
this method.

Constructing geometry on the basis of the deformation principle, it seems
that the formalism is absent, in general, as far as the theorems are absent.
(The geometers became accustomed that the geometric formalism appears only
in the form of theorem.) In reality the formalism appears in the implicit
form as a reference to the Euclidean geometry with its formalism. The
formalism of the Euclidian geometry is modified at the replacement of the
world function of the Euclidean space (Algorithms retain, only operand of
algorithms is changed). Besides, for a realization of such a modification
the Euclidean geometry is to be presented in the terms of the world
function. Mathematicians know slightly  this presentation, based on the
world function formalism. Despite the evidence of the deformation principle,
I obtained it, when the generalized geometry (T-geometry), based on its
application, has been already constructed. Furthermore, appearance of the
world function formalism (description of the Riemannian geometry in terms of
the world function) takes priority of the T-geometry construction.

\section{Application of T-geometry physical problems}

Besides, the overdetermination, appearing at a use of the conventional
(topological) method of the geometry construction, this methods is not
enough general. The method of the geometry construction based on the
deformation principle, admits one to construct the generalized geometry,
using deformation transforming the one-dimensional line into a surface.
Although the everyday experience does not provide an example of such
deformations, such deformations exist. Such a deformation may be described
as follows. Let us imagine a beam of straight thin elastic wires of the same
length. Being collected in a beam they represent a segment of the straight
line. Let us grasp the ends of wires by two hands and move hands one to
another. At such a deformation the wires diverge, conserving their lengths,
and form  an inswept surface. In other words, although deformations turning
a one-dimensional line into a surface are possible, the Riemannian geometry
and the metrical one ignore them, preferring to remain in the framework of
deformations, conserving one-dimensionality of the straight. A use of the
Euclidean concept of the straight, as a second base concept of the geometry,
and consideration of the one-dimensionality of the straight do not admit one
to prove out the deformation principle. The way to the deformation principle
lies through the construction of the mathematical formalism, based on the
application of the world function.

Applying the described above deformation to the construction of the
space-time geometry \cite{R2006}, one succeeded to construct such a
geometry, where the motion of free particles is primordially stochastic, and
intensity of the stochasticity depends on the particle mass. After proper
choice of the parameters of the space-time deformation (it depends on the
quantum constant) the statistical description of free stochastic particles
appears to be equivalent to the quantum description. \cite{R98}. It is
important, that the quantum principles are not used. In other words, the
quantum properties are descried by means of the correctly chosen space-time
geometry. It is impossible in the framework of the Riemannian geometry as
well as in the framework of any generalized geometry, constructed by means
of the conventional method of the generalized geometry construction.

\end{document}